\documentclass[11pt]{article}
\usepackage{amsmath, amsthm,amsfonts,amssymb,amscd}
\usepackage{color}
\usepackage{geometry}
\usepackage{tikz}
\usepackage{subfigure}
\numberwithin{equation}{section}
\usepackage[T1]{fontenc}
\usepackage{cite}
\usepackage{enumerate}

\usepackage{float}
\usepackage{soul}
\usepackage{marginnote}
\def\tto{\;{\lower 1pt \hbox{$\rightarrow$}}\kern -10pt
\hbox{\raise 2pt \hbox{$\rightarrow$}}\;}

\def\ra{\rangle}
\def\la{\langle}

\def\B{\mathbb B}

\def\h{\hfill\Box}
\def\R{\mathbb R}

\def\ox{\bar{x}}

\def\oy{\bar{y}}

\def\h{\hfill\triangle}

\def\ph{\varphi}

\newcounter{lk}

\begin{document}

\newtheorem{Theorem}{Theorem}[section]
\newtheorem{Proposition}[Theorem]{Proposition}
\newtheorem{Remark}[Theorem]{Remark}
\newtheorem{Lemma}[Theorem]{Lemma}
\newtheorem{Corollary}[Theorem]{Corollary}
\newtheorem{Definition}[Theorem]{Definition}
\newtheorem{Example}[Theorem]{Example}
\newtheorem{Fact}[Theorem]{Fact}
\renewcommand{\theequation}{\thesection.\arabic{equation}}
\normalsize
\def\proof{
\normalfont
\medskip
{\noindent\itshape Proof.\hspace*{6pt}\ignorespaces}}
\def\endproof{$\h$ \vspace*{0.1in}}

\title{\small \bf REGULAR CODERIVATIVE AND GRAPHICAL DERIVATIVE OF THE METRIC PROJECTION ONTO CLOSED BALLS IN HILBERT SPACES}
\date{}

\author{ Le Van Hien \footnote{Faculty of Pedagogy, Ha Tinh University, Ha Tinh city, Ha Tinh, Vietnam; email: hien.levan@htu.edu.vn}}

\maketitle

{\small \begin{abstract}
In this paper, we first establish a formula for exactly computing the regular coderivative  of the metric projection operator onto closed balls $r\B$ centered at the origin in Hilbert spaces. Then, this result is extended to metric projection operator onto any closed balls $\B(c,r)$, which has center $c$ in Hilbert space $H$ and with radius $r > 0$.  Finally, we give the formula for calculating the graphical derivative of the metric projection operator onto closed balls with center at arbitrarily given point in Hilbert spaces.
	
\end{abstract}}

\medskip

\normalsize
\section{Introduction}
\setcounter{equation}{0}

Let  $X$ be a general Banach space with norm $\| . \|$
and let $C$ be a nonempty subset of $X$. The metric projection  $P_C$
is defined by
$$P_C(x)=\{y\in C:\inf\limits_{z\in C}\|x-z\|=\|x-y\|\},\ \ \mbox{for all}\ \ x\in X.$$

The metric projection $P_C(x)$ is a set-valued mapping in the general case.  However, if $X$ is a Hilbert space and $C$ a nonempty closed convex subset, the metric projection operator is a well-defined single-valued mapping that has many useful properties, such as continuity, monotonicity, and non-expansiveness. Therefore, the metric projection operator becomes a useful and powerful tool in the optimization, computational mathematics, theory of equation, control theory, and others. Studying the properties of projection operators in Hilbert spaces is an interesting topic, attracting the attention of many people and many results in this direction have been established; see \cite{FP82, H77, K84, Li23, Li23(2), Li24, Li24(2), LLX23, M02, N95, S87, WY03} and the references therein.

Recently, investigating some characteristic properties, the Gâteaux directional differentiability, the Fréchet differentiability and the strict  Fréchet differentiabilityof projection operators has received special attention from Li et al. In 2024, Khan and Li \cite{KL24} studied some characteristic properties of the projection operator onto a specific subspace of the Banach space. In \cite{Li23, LLX23}, Li et al. studied the Gâteaux directional differentiability and proved some properties of  Gâteaux directional derivatives of $P_C$ in both Banach spaces and Hilbert spaces. In particular, when $C$ are the closed and convex subsets, or the closed balls, or the closed and convex cones (including proper closed subspaces) in Hilbert spaces, the exact Gâteaux directional derivatives of $P_C$ are provided in \cite{Li23,LLX23}. Li \cite{Li24} proved strict Fréchet differentiability of the metric projection operator $P_{r\B}$ onto closed balls $r\B$ centered at the origin in Hilbert spaces and onto positive cones in Euclidean spaces. In this way, Hien \cite{H24} studied the strict Fréchet differentiability of the metric projection $P_{\B(c,r)}$, which is onto closed balls with center at arbitrarily given point $c$ in the spaces and the metric projection operator onto the second-order cones in Euclidean spaces.

The regular coderivative and the graphical derivative are powerful tool in variational analysis and its applications. These tools can be used to investigate the stability and sensitivity ofconstraint and variational systems and, more generally, generalized equations \cite{CHN18, GM15, MOR15, OS08} or, even, characterize some nice properties of set-valued mappings  \cite{DQZ06, DR14, M18, RW98}. Although the regular coderivative and the graphical derivative are an important key to solving many problems, calculating it is not a simple thing, attracting the attention of many people and many interesting results have been established; see \cite{CH17, GM15, M18, HOS12, MOR15, RW98, YZ17} and the references therein.

In this paper, we are interested in the regular coderivative and the graphical derivative of the metric projection operator $P_C$, where $C$ is the closed ball in Hilbert spaces. For $C$  is the closed balls $r\B$ centered at the origin in Hilbert spaces, Li \cite{Li23(2)} has given some properties on the regular coderivative of $P_{r\B}$. These results are then also generalized to the metric projection onto closed balls centered at the origin in uniformly convex and uniformly smooth Banach spaces \cite{Li24(2)}. However, in these works, the results of calculating the regular coderivative of $P_{r\B}$ are only given in some special cases and there is still no information about the regular coderivative of the metric projection operator $P_{\B(c,r)}$.


The aim of this paper is to establish a formula for calculating the regular coderivative and the graphical derivative of the metric projection operator onto closed balls in Hilebert spaces. More precisely, we first establish aformula for exactly computing the regular coderivative  of the metric projection operator onto closed balls $r\B$ centered at the origin in Hilbert spaces. Then, this result is extended to metric projection $P_{\B(c,r)}$, which is onto closed balls with center at arbitrarily given point $c$ in the spaces.  Finally, we give the formula for calculating the graphical derivative of the metric projection operator onto closed balls in Hilbert spaces.

The structure of the paper is as follows.   In section 2,  we recall some preliminary materials. In section 3, we present results on the regular coderivative of the metric projection operator onto closed balls in Hilbert spaces. Then, in section 4, we give the exact calculation formula for the graphical derivative of the metric projection operator onto closed balls in Hilbert spaces. Finally, we conclude the paper in section 5 where we discuss some perspectives of the obtained results and future works.

\section{Preliminaries}
\setcounter{equation}{0}
 
  We first give the following notations that will be used throughout the paper. Let $X$ be Banach spaces, we denote by $\|.\|$ the norm in $X$. Given a set $C \subset X,$ we denote by $C^o$, $\mbox{bd}C$ its interior, boundary respectively. The closed ball with center $\ox$ and radius $r>0$ is denoted by $\B(\ox,r)$.   Define the identity mapping on $X$ by $I_X:X\to X$ with $I_X(x)=x,\ \mbox{for all}\ x\in X.$  Denote $t\downarrow 0$ means that $t\to 0$ with $t>0$ and $u\xrightarrow{C}\ox$   means that $u\rightarrow \ox$ with $u\in C$.


Below are basic notions and facts from variational analysis, which are frequently used in the following; see \cite{A99, Li24, M18, RW98, S96} for more details.
\begin{Definition}{\rm(see \cite{A99, Li24, M18, RW98, S96})
Let $f:X\to Y$ be a single-valued map between two Banach spaces $X, Y$ and $ \Omega\subset X$. 

$(i)$ If there is a  real constant $L\geq 0$ such that $$\|f(x_1)-f(x_2)\|\leq L\|x_1-x_2\|,\ \ \mbox{ for all}\ x_1, x_2 \in \Omega,$$
then $f$ is said to {\it Lipschitz continuous} on $\Omega$.

$(ii)$ The mapping $f$ is called {\it Gâteaux directionally differentiable} at point $x\in X$ along direction $w\in X$ if the following limit exists $$\lim\limits_{t\downarrow0}\dfrac{f(x+tw)-f(x)}{t}:=f'(x)(w).$$
Then, $f'(x)(w)$ is said to be the {\it Gâteaux directional derivative} of $f$ at point $x$ along direction $w$.

The mapping $f$ is called {\it Gâteaux directionally differentiable on $\Omega$} if $f$ is Gâteaux directionally differentiable at every point $x\in \Omega$. 

$(iii)$ The mapping $f$ is called {\it Fréchet differentiable} at $\ox \in X$ if there is a linear continuous mapping $\nabla f(\ox):X\to Y$ such that
$$\lim\limits_{h\to 0}\dfrac{\|f(\ox+h)-f(\ox)-\nabla f(\ox)(h)\|}{\|h\|}=0.$$
Then, $\nabla f(\ox)$  is said to be the {\it Fréchet
	derivative} of $f$ at $\ox$.
	
In particular, the mapping $f$ is called {\it strictly Fréchet differentiable} at $\ox \in X$ if
$$\lim\limits_{(u,v)\to (\ox,\ox)}\dfrac{\|f(u)-f(v)-\nabla f(\ox)(u-v)\|}{\|u-v\|}=0.$$

The mapping $f$ is called {\it Fréchet differentiable (strictly Fréchet differentiable) on $\Omega$} if $f$ is Fréchet differentiable (strictly Fréchet differentiable, respectively) at every point $x\in \Omega$.	
	
}
\end{Definition}




\begin{Remark} {\rm 
		The strict Fréchet differentiability entails the Fréchet differentiability; the Fréchet differentiability entails the Gâteaux directional differentiability and for any $\ox, u\in ~X$, we have $f'(\ox)(w)=\nabla f(\ox)(w)$. However, the converse is not true (see \cite{H24, Li24}).
}\end{Remark}
\begin{Definition} {\rm (see \cite{M18,RW98})
	Let $\Omega$ be a nonempty subset of  a Hilbert space $H$. The {\it (Fréchet) regular normal cone} to $\Omega$ at $\ox\in\Omega$ is the set $\widehat N_\Omega(\ox)$ given by
	$$\widehat N_\Omega(\ox):=\Big\{z\in H| \limsup\limits_{u\xrightarrow{\Omega} \ox}\dfrac{\la z,u-\ox\ra}{\|u-\ox\|} \leq 0 \Big\}.$$}\end{Definition}
	If $\ox\not\in\Omega$, one puts
	$\widehat N_\Omega(\ox)=\emptyset$ by convention. When the set $\Omega$ is convex, the regular normal cone reduce to the normal cone in the sense of classical convex analysis.
	
	\begin{Definition} {\rm (see \cite{M18,RW98})
	Let $H$ be Hilbert space and $\Phi: H\rightrightarrows H$ be a set-valued mapping with its graph {\rm gph}$\Phi: =\{(x,y)|y\in\Phi(x)\}$ and its domain {\rm Dom}$\Phi:=\{x|\Phi(x)\not=\emptyset\}$. The {\it (Fréchet) regular coderivative} of $\Phi$ at a given point $(\ox,\oy)\in H\times H$ is the set-valued mapping $\widehat D^\ast\Phi(\ox,\oy):H\rightrightarrows H$ defined by
$$\widehat D^\ast\Phi(\ox,\oy)(w):=\big\{v\in H| (v,-w)\in \widehat N_{\mbox{\rm gph}\Phi}(\ox, \oy)\}\ \mbox{for all}\ w\in H.$$ }
	\end{Definition}
	In the case $\Phi(\ox)=\{\oy\},$ one writes  $\widehat D^\ast\Phi(\ox)$ for  $\widehat D^\ast\Phi(\ox,\oy)$.\\
Note that {\rm gph}$\Phi \subset H\times H$, where $H\times H$ be the orthogonal product of $H$ equipped with the inner product $\la. , .\ra_{H\times H}$, which is abbreviated as $\la. , .\ra$.  It is defined by, $$\la(x,y), (u,v)\ra=\la x,v\ra +\la y,v\ra,\ \ \mbox{for any}\  (x,y), (u,v)\in H \times H.$$
Let $\| . \|_{H\times H}$ be the norm in $H\times H$ induced by the inner product $\la. , .\ra$. Then, we have
$$\| (x,y) \|_{H\times H}=\sqrt{\la(x,y), (x,y)\ra} =\sqrt{\|x\|^2+\|y\|^2}\ \ \mbox{for all}\ x,y\in H.$$

	The following lemma provides a useful formula to calculate the regular coderivatives of single-valued mappings on Hilbert spaces.	
	
\begin{Lemma} \label{Lem1}
	Let $H$ be a Hilbert space, $f: H\rightarrow H$ be a Lipschitz  continuous mapping on $H$. Then, the regular coderivatives of $f$ at a point $\ox\in H$ satisfies that, for any $y\in H$,
	$$\widehat D^\ast f(\ox)(y)=\left\{z\in H\big|\limsup\limits_{u\rightarrow \ox}\dfrac{\la z,u-\ox\ra-\la y,f(u)-f(\ox)\ra}{\|u-\ox\|} \leq 0 \right\}.$$
\end{Lemma}
{\bf Proof.} Thanks to the continuity of $f$, by the definition of the regular coderivative, we have
$$\begin{array}{rl}z\in \widehat D^\ast f(\ox)(y) &\Longleftrightarrow (z,-y)\in\widehat N_{\mbox{\rm gph}f}(\ox)\\
	&\Longleftrightarrow \limsup\limits_{(u,v)\xrightarrow{\mbox{\rm gph}f} (\ox, f(\ox))}\dfrac{\la (z,-y), (u,v)-(\ox, f(\ox))\ra}{\|(u,v)-(\ox, f(\ox))\|_{H\times H}} \leq 0\\
		&\Longleftrightarrow \limsup\limits_{u\rightarrow\ox,
		v=f(u)}\dfrac{\la z, u-\ox\ra -\la y, v- f(\ox))\ra}{\|(u-\ox,v- f(\ox))\|_{H\times H}} \leq 0\\
		&\Longleftrightarrow \limsup\limits_{u\rightarrow\ox}\dfrac{\la z, u-\ox\ra -\la y, f(u)- f(\ox))\ra}{\sqrt{\|u-\ox\|^2+\|f(u)-f(\ox)\|^2}} \leq 0.
\end{array}$$
Notice that $f$ be a Lipschitz  continuous mapping on $H$, there exists $L\geq0$ such that $\|f(x_1)-f(x_2)\|\leq L\|x_1-x_2\|,\ \ \mbox{ for all}\ x_1, x_2 \in H.$ 

So,
$$
\|u-\ox\|\leq\sqrt{\|u-\ox\|^2+\|f(u)-f(\ox)\|^2}\leq \sqrt{1+L^2}\|u-\ox\|
$$
for all $\ox, u\in H.$\\
We consider the following two cases:\\
{\it Case 1.} If $\la z, u-\ox\ra -\la y, f(u)- f(\ox))\ra\geq 0$ then
$$\begin{array} {rl}\dfrac{\la z, u-\ox\ra -\la y, f(u)- f(\ox))\ra}{\sqrt{1+L^2}\|u-\ox\|}&\leq\dfrac{\la z, u-\ox\ra -\la y, f(u)- f(\ox))\ra}{\sqrt{\|u-\ox\|^2+\|f(u)-f(\ox)\|^2}}\\
	&\leq\dfrac{\la z, u-\ox\ra -\la y, f(u)- f(\ox))\ra}{\|u-\ox\|}.\end{array}$$
{\it Case 2.} If $\la z, u-\ox\ra -\la y, f(u)- f(\ox))\ra< 0$ then
$$\begin{array} {rl}\dfrac{\la z, u-\ox\ra -\la y, f(u)- f(\ox))\ra}{\|u-\ox\|}&\leq\dfrac{\la z, u-\ox\ra -\la y, f(u)- f(\ox))\ra}{\sqrt{\|u-\ox\|^2+\|f(u)-f(\ox)\|^2}}\\
	&\leq\dfrac{\la z, u-\ox\ra -\la y, f(u)- f(\ox))\ra}{\sqrt{1+L^2}\|u-\ox\|}.\end{array}$$
Therefore, $$\begin{array}{rl}z\in \widehat D^\ast f(\ox)(y) 
	&\Longleftrightarrow \limsup\limits_{u\rightarrow\ox}\dfrac{\la z, u-\ox\ra -\la y, f(u)- f(\ox))\ra}{\sqrt{\|u-\ox\|^2+\|f(u)-f(\ox)\|^2}} \leq 0\\
	&\Longleftrightarrow \limsup\limits_{u\rightarrow\ox}\dfrac{\la z, u-\ox\ra -\la y, f(u)- f(\ox))\ra}{\|u-\ox\|} \leq 0.
\end{array}$$\hfill $\square$

The following results, which were presented in \cite{M18}, are useful for us to obtain the main result in the next section.

\begin{Proposition}\label{Pro13} {\rm (see \cite[Theorem 1.38]{M18})}
	Let $H$ be Hilbert space, $f: H\to H$ be Fréchet differentiable  at $\ox\in H.$   Then
	$$\widehat D^\ast f(\ox)(y)=\big\{\nabla f(\ox)(y)\big\}\ \ \mbox{for all}\ y\in H.$$
\end{Proposition}
\begin{Proposition}\label{Pro23} {\rm (see \cite[Theorem 1.62]{M18})}
	Let $H$ be Hilbert space, $f: H\to H$ be Fréchet differentiable  at $\ox\in H,$ and let $g: H\to H$ be an arbitrary mapping.  Then
	$$\widehat D^\ast(f + g)(\ox)(y)=\nabla f(\ox)+\widehat D^\ast g(\ox).$$
\end{Proposition}
\begin{Proposition}\label{Pro24} {\rm (see \cite[Theorem 1.66]{M18})}
Let $H$ be Hilbert space, $f, g: H\to H$. Assume that $f$ is strictly Fréchet differentiable at $\ox$ with the surjective derivative $\nabla f(\ox)$. Then, we have
	\begin{equation}\label{Pro26} \widehat D^\ast (g\circ f)(\ox)=\nabla f(\ox)\widehat D^\ast g(f(\ox)).\end{equation}
\end{Proposition}

 \section{Regular coderivative of the metric projection onto balls in Hilbert spaces}\label{Sec3}
\setcounter{equation}{0}

In this section, let $(H,\| . \|)$ be a real Hilbert space with inner product $\la\ .\ \ra$. For any $c\in H$ and $r>0$, let $\B(c,r)$ denote the closed ball in $H$ with  centered at $c$ and  radius $r$. Notation $\B(\theta,1)=\B$ and $\B(\theta,r)=r\B.$
Let $P_{\B(c,r)}:H\to \B(c,r)$ be the  metric projection operator onto  $\B(c,r)$ in $H$. 

It is well-known that $P_{\B(c,r)}$ is a well-defined single-valued mapping, and
$$P_{\B(c,r)}(x)=\begin{cases}
	x\ \ \ \ \ \ \ \ \ \ \ \ \ \ \ \ \ \ \ \ \ \ \ \ \  \mbox{if} \ x\in \B(c,r);\\
	c+\dfrac{r}{\|x-c\|}(x-c)\ \ \ \mbox{if} \ x\in H\backslash\B(c,r).
\end{cases}$$
In particular, 
$$P_{r\B}(x)=\begin{cases}
	x\ \ \ \ \ \ \ \ \ \ \ \ \  \mbox{if} \ x\in r\B;\\
	\dfrac{rx}{\|x\|}\ \ \ \ \ \ \ \  \ \mbox{if} \ x\in H\backslash r\B.
\end{cases}$$
Furthermore, $P_{\B(c,r)}$ is a  Lipschitz continuous mapping on $H$ with Lipschitz constant L=1, i.e 
$$\|P_{\B(c,r)}(x_1)-P_{\B(c,r)}(x_2)\|\leq \|x_1-x_2\|,\ \mbox{for all}\ \ x_1,x_2\in H.$$
For any $x\in{\rm\mbox{bd}}\B(c,r)$, two subsets $x^\uparrow_{(c,r)}$ and $x^\downarrow_{(c,r)}$ of $H\backslash\{\theta\}$ are defined by: for $v\in H$ with $v\not = \theta$,
$$x^\uparrow_{(c,r)}=\{v\in H\backslash\{\theta\}|\, \mbox{there is}\  \delta>0 \ \mbox{such that}\ \|(x+tv)-c\|\geq r,\ \mbox{for all}\ t\in (0,\delta)\}.$$
$$x^\downarrow_{(c,r)}=\{v\in H\backslash\{\theta\}|\, \mbox{there is}\  \delta>0 \ \mbox{such that}\ \|(x+tv)-c\|<r,\ \mbox{for all}\ t\in (0,\delta)\}.$$
In particular, we write
$$x^\uparrow:=x^\uparrow_{(\theta,1)}\ \  \mbox{and}\ \ x^\downarrow:=x^\downarrow_{(\theta,1)}.$$
By \cite[Lemma 5.1]{Li23}, we have $x^\uparrow_{(c,r)}\cup x^\downarrow_{(c,r)} =\emptyset$ and $x^\uparrow_{(c,r)}\cap x^\downarrow_{(c,r)} =H\backslash\{\theta\}$, for all $x\in{\rm\mbox{bd}}\B(c,r)$ with $c\in H$ and $r>0$.

For any $x\in H\backslash \{\theta\}$, let $S(x)$ be the one-dimensional subspace of $H$ generated by $x$. Let $O(x) $ denote the orthogonal subspace of $x$ in $H$. Then, $H$ has the following orthogonal decomposition
$$H=S(x)\oplus O(x).$$ 
Specifically, $x\in H\backslash \{\theta\}$ and for all $u\in H$, we have the following orthogonal representation
$$u=\dfrac{\la x, u\ra}{\|x\|^2}x+\left(u-\dfrac{\la x, u\ra}{\|x\|^2}x\right).$$
Therefore, for any fixed $x\in H\backslash\{\theta\}$, we define a function $a(x, . ): H\to \R$ by
$$a(x, u)=\dfrac{\la x, u\ra}{\|x\|^2}, \ \ \ \ \mbox{for all}\ u\in H.$$
And a mapping  $o(x, . ):H\to O(x)$ by
$$o(x, u)=u-\dfrac{\la x, u\ra}{\|x\|^2}x, \ \ \ \ \mbox{for all}\ u\in H.$$
We see that \\

 $\blacklozenge$\ \ \ \ \  $u=a(x,u)x+o(x,u)\ \ \ \mbox{and}\ \ \ \big\la a(x,u)x, o(x,u)\big\ra=0, \ \ \mbox{for all}\ u\in H.$

$\blacklozenge$\ \ \ \ \ $a(x, . ), o(x, . )$ are continuous linear mappings, and $o(x,x)=\theta$, $o(x,x+h)=o(x,h)$ for all $h\in H.$ 

$\blacklozenge$\ \ \ \ \  $\la u,v\ra=a(x,u)a(x,v)\|x\|^2+\la o(x,u), o(x,v)\ra, \ \mbox{for all}\ u,v\in H.$

$\blacklozenge$\ \ \ \ \  $\|u\|^2=(a(x,u))^2\|x\|^2+\|o(x,u)\|^2, \ \mbox{for all}\ u\in H.$

$\blacklozenge$\ \ \ \ \  $\|u+v\|^2=\big((a(x,u))+a(x,v)\big)^2\|x\|^2+\|o(x,u)+o(x,v)\|^2, \ \mbox{for all}\ u, v\in H.$

$\blacklozenge$\ \ \ \ \ $a(x,u)\to 1$ and $o(x,u)\to \theta$ as $u\to x,\ \ \mbox{for}\ u\in H.$

For the sake of simplicity, $a(x,u)$ and $o(x,u)$ are abbreviated as $a(u)$ and $o(u)$, respectively.\\

The following result provides the formula for calculating the regular coderivatives of the metric projection onto a closed ball with centered at the origin.

	\begin{Theorem} \label{Thm1} 
		Let $H$ be a Hilbert space. For any  $r>0$, the regular coderivative of $P_{r\B}:H\rightarrow r\B$   at a point $\ox\in H$ is given below.

		(i)  If $\ox\in r\B^o$, $$\widehat D^\ast P_{r\B}(\ox)(y)=\{y\} \ \ \mbox{for every}\ y\in H.$$
		
		(ii) If $\ox\in H\backslash r\B$, then
		$$\widehat D^\ast P_{r\B}(\ox)(y)=\left\{\dfrac{r}{\|\ox\|}\Big(y-\dfrac{\la \ox, y\ra}{\|\ox\|^2}\ox\Big)\right\}\ \ \mbox{for every} \ y\in H.$$

		(iii) If $\ox\in \mbox{\rm bd}r\B$, then
		
\begin{equation}\label{kq}\widehat D^\ast P_{r\B}(\ox)(y)=\left\{z\in H: \la y,\ox\ra\leq \la z,\ox\ra \leq 0, \  y-z=\dfrac{\la y-z,\ox\ra}{\|\ox\|^2}\ox\right\}.\end{equation}

	\end{Theorem}
	
	{\bf Proof.} Parts $(i, ii)$ follow from \cite[Theorem 3.3]{Li24} and Proposition \ref{Pro13}. So, we only prove part $(iii)$.
	
  Let $y\in H$, thanks to the Lipschitz continuity of $P_{r\B}$, by Lemma \ref{Lem1}, we have
\begin{equation}\label{39}z\in \widehat D^\ast P_{r\B}(\ox)(y) \iff \limsup\limits_{u\rightarrow \ox}\dfrac{\la z,u-\ox\ra-\la y,P_{r\B}(u)-P_{r\B}(\ox)\ra}{\|u-\ox\|} \leq 0.\end{equation}

Suppose that $z\in \widehat D^\ast P_{r\B}(\ox)(y)$. We take a directional line segment in the limit in \eqref{39}, $u = (1-\delta)\ox$ for $\delta\downarrow 0$ with $\delta <1$. It follows that
$$\begin{array}{rl}
	0&\geq \limsup\limits_{\delta\downarrow 0}\dfrac{\la z,-\delta\ox\ra-\la y,-\delta\ox\ra}{\delta\|\ox\|}\\
	&=\limsup\limits_{\delta\downarrow 0}\dfrac{\la y-z,\ox\ra}{\|\ox\|}\\
	&=\limsup\limits_{\delta\downarrow 0}\dfrac{\la (a(y)-a(z))\ox,\ox\ra}{\|\ox\|}\\
	&= (a(y)-a(z))r.
	\end{array}$$
This implies that $a(y)\leq a(z)$, it means\begin{equation}\label{10} \la y,\ox\ra\leq \la z,\ox\ra \end{equation}

We take a directional line segment in the limit in \eqref{39}, $u = (1+\delta)\ox$ for $\delta\downarrow 0$ with $\delta <1$. It follows that
$$\begin{array}{rl}
	0&\geq \limsup\limits_{\delta\downarrow 0}\dfrac{\la z,\delta\ox\ra-\la y,\ox-\ox\ra}{\delta\|\ox\|}\\
	&=\limsup\limits_{\delta\downarrow 0}\dfrac{\la a(z)\ox+o(z),\ox\ra}{\|\ox\|}\\
	&=\limsup\limits_{\delta\downarrow 0}\dfrac{\la a(z)\|\ox\|^2}{\|\ox\|}\\
	&= a(z)r.
\end{array}$$
This implies that $a(z)\leq 0$, it means  \begin{equation}\label{11} \la z,\ox\ra\leq 0. \end{equation}

Assume, by the way of contradiction, that \begin{equation} \label{311}
o(z)\not=o(y).
\end{equation}

Under the hypothesis \eqref{311}, we take a directional line segment in the limit in \eqref{39}, $u=\ox-\delta\big(o(y)-oz\big)$, for $\delta\downarrow 0$ with $0 < \delta < 1$. We consider the following two cases:\\
{\it Case 1.} $o(z)-o(y)\in\ox^\uparrow.$ We have $P_{r\mathbb{B}}(u)=\dfrac{r}{\|u\|}u$ and
$$\begin{array}{rl}
	0&\geq \limsup\limits_{\delta\downarrow 0}\dfrac{\la z,u-\ox\ra-\la y,P_{r\mathbb{B}}(\ox-\delta\big(o(y)-o(z)\big))-P_{r\mathbb{B}}(\ox)\ra}{\delta\|o(z)-o(y)\|}\\
	&=\limsup\limits_{\delta\downarrow 0}\dfrac{\delta\big\la a(z)\ox+o(z),o(z)-o(y)\big\ra-\big\la a(y)\ox+o(y),\dfrac{r}{\|u\|}(\ox-\delta\big(o(y)-o(z)\big))-\ox\big\ra}{\delta\|o(z)-o(y)\|}\\
	&=\limsup\limits_{\delta\downarrow 0}\dfrac{\delta\big\la o(z),o(z)-o(y)\big\ra- a(y)\|\ox\|^2\bigg(\dfrac{r}{\|u\|}-1\bigg)+\dfrac{r\delta}{\|u\|}\big\la o(y), o(z)-o(y)\big\ra}{\delta\|o(z)-o(y)\|}\\
	&=\limsup\limits_{\delta\downarrow 0}\dfrac{\delta\big\la o(z),o(z)-o(y)\big\ra- a(y)r^2\dfrac{r^2-\|u\|^2}{\|u\|(r+\|u\|)}+\dfrac{r\delta}{\|u\|}\big\la o(y), o(z)-o(y)\big\ra}{\delta\|o(z)-o(y)\|}\\		&=\limsup\limits_{\delta\downarrow 0}\dfrac{\delta\big\la o(z),o(z)-o(y)\big\ra- a(y)r^2\dfrac{r^2-\|\ox\|^2-\delta^2\|o(y)-o(z)\|^2}{\|u\|(r+\|u\|)}+\dfrac{r\delta}{\|u\|}\big\la o(y), o(z)-o(y)\big\ra}{\delta\|o(z)-o(y)\|}\\
&=\limsup\limits_{\delta\downarrow 0}\dfrac{\delta\big\la o(z),o(z)-o(y)\big\ra+ a(y)r^2\delta^2\dfrac{\|o(y)-o(z)\|^2}{\|u\|(r+\|u\|)}+\dfrac{r\delta}{\|u\|}\big\la o(y), o(z)-o(y)\big\ra}{\delta\|o(z)-o(y)\|}\\
&=\limsup\limits_{\delta\downarrow 0}\dfrac{\big\la o(z),o(z)-o(y)\big\ra+ a(y)r^2\delta\dfrac{\|o(y)-o(z)\|^2}{\|u\|(r+\|u\|)}+\dfrac{r}{\|u\|}\big\la o(y), o(z)-o(y)\big\ra}{\|o(z)-o(y)\|}\\	
&=\dfrac{\big\la o(z),o(z)-o(y)\big\ra+ \dfrac{r}{\|\ox\|}\big\la o(y), oz-oy\big\ra}{\|o(z)-o(y)\|}\\
&=\dfrac{\big\la o(z)-oy,o(z)-o(y)\big\ra}{\|o(z)-o(y)\|}\\
&=\|o(z)-o(y)\|.
\end{array}$$

{\it Case 2.} $o(z)-o(y)\in\ox^\downarrow.$ We have $P_{r\mathbb{B}}(u)=u$ and
$$\begin{array}{rl}
	0&\geq \limsup\limits_{\delta\downarrow 0}\dfrac{\la z,\ox-\delta(o(y)-o(z))-\ox\ra-\la y,P_{r\mathbb{B}}(\ox-\delta(o(y)-o(z)))-P_{r\mathbb{B}}(\ox)\ra}{\|\ox-\delta(o(y)-o(z))-\ox\|}\\

	&= \limsup\limits_{\delta\downarrow 0}\dfrac{\la z,\delta(o(z)-o(y))\ra-\la y, \ox-\delta\big(o(y)-o(z)\big)-\ox\ra}{\delta\|o(y)-o(z)\|}\\
	&= \limsup\limits_{\delta\downarrow 0}\dfrac{\la z-y,\delta(o(z)-o(y))\ra}{\delta\|o(y)-o(z)\|}\\
	&=\dfrac{\la z-y,o(z)-o(y)\ra}{\|o(y)-o(z)\|}\\
	&=\dfrac{\la o(z)-o(y),o(z)-o(y)\ra}{\|o(y)-o(z)\|}\\
		&=\|o(y)-o(z)\|.
\end{array}
$$

In both cases we obtain a contradiction with \eqref{311}, so $o(z)=o(y)$, it means \begin{equation} \label{12}
 y-z=\dfrac{\la y-z,\ox\ra}{\|\ox\|^2}\ox.
\end{equation}
By \eqref{10}, \eqref{11} and \eqref{12} together, we have 
\begin{equation} \label{c}
	z\in \widehat D^\ast P_{r\B}(\ox)(y)\Longrightarrow  \la y,\ox\ra\leq \la z,\ox\ra \leq 0\ \  \mbox{and}\ \  y-z=\dfrac{\la y-z,\ox\ra}{\|\ox\|^2}\ox.
\end{equation}

Next, we prove the converse, assuming $\la y,\ox\ra\leq \la z,\ox\ra \leq 0\ \  \mbox{and}\ \  y-z=\dfrac{\la y-z,\ox\ra}{\|\ox\|^2}\ox$, it means $a(y)\leq a(z)\leq 0\ \  \mbox{and}\ \ o(z)=o(y)$. We need to show that
$$\limsup\limits_{u\rightarrow \ox}\dfrac{\la z,u-\ox\ra-\la y,P_{r\B}(u)-P_{r\B}(\ox)\ra}{\|u-\ox\|}\leq 0.$$
We consider the following cases:

{\it Case 1.} $u\in r\mathbb{B}$. We have  $P_{r\mathbb{B}}(u)=u, \ a(u)<1$ and 
\begin{equation} \label{1}\begin{array}{rl}\limsup\limits_{u\xrightarrow {r\mathbb{B}}\ox}\dfrac{\la z,u-\ox\ra-\la y,P_{r\B}(u)-P_{r\B}(\ox)\ra}{\|u-\ox\|}&=\limsup\limits_{u\xrightarrow {r\mathbb{B}}\ox}\dfrac{\la z,u-\ox\ra-\la y,u-\ox\ra}{\|u-\ox\|}\\
	&=\limsup\limits_{u\xrightarrow {r\mathbb{B}}\ox}\dfrac{\la z-y,u-\ox\ra}{\|u-\ox\|}\\
	&=\limsup\limits_{u\xrightarrow {r\mathbb{B}}\ox}\dfrac{\big\la \big(a(z)-a(y)\big)\ox,a(u)\ox+o(u)-\ox\big\ra}{\|u-\ox\|}\\
		&=\limsup\limits_{u\xrightarrow {r\mathbb{B}}\ox}\dfrac{ \big(a(z)-a(y)\big)\big(a(u)-1\big)\|\ox\|^2}{\|u-\ox\|}\\
		&\leq 0.
	\end{array}\end{equation}
	{\it Case 2.} $u\in H\backslash r\mathbb{B}$. We have  $P_{r\mathbb{B}}(u)=\dfrac{r}{\|u\|}u$ and 
	\begin{equation} \label{3.7}
	\begin{array}{rl}&\limsup\limits_{u\xrightarrow {H\backslash r\mathbb{B}}\ox}\dfrac{\la z,u-\ox\ra-\la y,P_{r\B}(u)-P_{r\B}(\ox)\ra}{\|u-\ox\|}\\
		&=\limsup\limits_{u\xrightarrow {H\backslash r\mathbb{B}}\ox}\dfrac{\la z,u-\ox\ra-\Big\la y,\dfrac{r}{\|u\|}u-\ox\Big\ra}{\|u-\ox\|}\\
		&=\limsup\limits_{u\xrightarrow {H\backslash r\mathbb{B}}\ox}\dfrac{\la a(z)\ox+o(z),a(u)\ox+o(u)-\ox\ra-\Big\la a(y)\ox+o(y),\dfrac{r}{\|u\|}\big(a(u)\ox+o(u)\big)-\ox\Big\ra}{\|u-\ox\|}\\
	&=\limsup\limits_{u\xrightarrow {H\backslash r\mathbb{B}}\ox}\dfrac{ a(z)\big(a(u)-1\big)\|\ox\|^2+\big\la o(z),o(u)\big\ra-a(y)\Big(\dfrac{r a(u)}{\|u\|}-1\Big)\|\ox\|^2-\dfrac{r}{\|u\|}\big\la o(y), o(u)\big\ra}{\|u-\ox\|}\\
	&=\limsup\limits_{u\xrightarrow {H\backslash r\mathbb{B}}\ox}\dfrac{ a(z)\big(a(u)-1\big)r^2-a(y)\Big(\dfrac{r a(u)}{\|u\|}-1\Big)r^2+\Big(1-\dfrac{r}{\|u\|}\Big)\big\la o(y), o(u)\big\ra}{\|u-\ox\|}\\
	\end{array}\end{equation}

Notice that \begin{equation}\label{3.8}\begin{array}{rl}0\leq \Bigg\|\dfrac{ \Big(1-\dfrac{r}{\|u\|}\Big)\big\la o(y), o(u)\big\ra}{\|u-\ox\|}\Bigg\|&=\dfrac{ \Big(1-\dfrac{r}{\|u\|}\Big)\Big|\big\la o(y), o(u)\big\ra\Big|}{\big\|\big(a(u)-1\big)\ox+o(u)\big\|}\\
	&\leq \dfrac{ \Big(1-\dfrac{r}{\|u\|}\Big)\Big|\big\la o(y), o(u)\big\ra\Big|}{\|o(u)\|}\\
		&\leq \dfrac{ \Big(1-\dfrac{r}{\|u\|}\Big)\big\|o(y)\big\| \big\|o(u)\big\|}{\|o(u)\|}\\
			&=  \Big(1-\dfrac{r}{\|u\|}\Big)\big\|o(y)\big\| \\
& \rightarrow 0\ \ \ \ \ \ \ \ \ \ \ \ \ \ \  \ \ \ \ \ \ \ \ \ \ \ \ \ \ \ \ \ \ \ \ \mbox{as}\ u\rightarrow\ox,
\end{array}
\end{equation}
and \begin{equation} \label{3.9}
a(z)\big(a(u)-1\big)r^2-a(y)\Big(\dfrac{r a(u)}{\|u\|}-1\Big)r^2\leq 0.
\end{equation}
Indeed, if $a(u)<1$, then $\dfrac{r a(u)}{\|u\|}-1\leq a(u)-1<0.$ Thus, $a(z)\big(a(u)-1\big)\leq a(y)\Big(\dfrac{r a(u)}{\|u\|}-1\Big)$ and $a(z)\big(a(u)-1\big)r^2-a(y)\Big(\dfrac{r a(u)}{\|u\|}-1\Big)r^2\leq 0$. If $a(u)>1$, then $a(z)\big(a(u)-1\big)r^2\leq 0$. We have $$\begin{array} {rl}
	a(y)\Big(\dfrac{r a(u)}{\|u\|}-1\Big)r^2&=	a(y)r^2\dfrac{r a(u)-\|u\|}{\|u\|}\\
	&=a(y)r^2\dfrac{r^2 a(u)^2-\|u\|^2}{\|u\|\big(r a(u)+\|u\|\big)}\\
	&=a(y)r^2\dfrac{-\|o(u)\|^2}{\|u\|\big(r a(u)+\|u\|\big)}\\
	&\geq 0.
\end{array}
$$
So, $a(z)\big(a(u)-1\big)r^2-a(y)\Big(\dfrac{r a(u)}{\|u\|}-1\Big)r^2\leq 0$. \\

Combining \eqref{3.7} with \eqref{3.8} and \eqref{3.9}, we obtain 
\begin{equation} \label{2}\limsup\limits_{u\xrightarrow {H\backslash r\mathbb{B}}\ox}\dfrac{\la z,u-\ox\ra-\la y,P_{r\B}(u)-P_{r\B}(\ox)\ra}{\|u-\ox\|}\leq 0.\end{equation}
By \eqref{1} and \eqref{2} together, we have
$$\limsup\limits_{u\rightarrow \ox}\dfrac{\la z,u-\ox\ra-\la y,P_{r\B}(u)-P_{r\B}(\ox)\ra}{\|u-\ox\|}\leq 0.$$
This means \begin{equation} \label{d}
\la y,\ox\ra\leq \la z,\ox\ra \leq 0\ \  \mbox{and}\ \  y-z=\dfrac{\la y-z,\ox\ra}{\|\ox\|^2}\ox\Longrightarrow 	z\in \widehat D^\ast P_{r\B}(\ox)(y).
\end{equation}
Then, by \eqref{c} and \eqref{d}, \eqref{kq} is proved.\hfill $\square$
\\

Using Theorem \ref{Thm1}, we obtain the following result, which was given by  Li in \cite{Li23(2)}.
\begin{Corollary} {\rm (see \cite[Theorem 3.2]{Li23(2)}) }
Let $H$ be a Hilbert space. For any $r > 0$, if $\ox\in \mbox{\rm bd}r\mathbb{B}$, then

(a) $ \widehat D^\ast P_{r\B}(\ox)(\theta)=\{\theta\};$\\

(b) For any $y\in H\backslash\{\theta\}$, we have
$$\theta\in \widehat D^\ast P_{r\B}(\ox)(y)\Longleftrightarrow y=\dfrac{\la y,\ox\ra}{\|\ox\|^2}\ox\ \ \mbox{with}\ \la y,\ox\ra\leq 0;$$

(c) $\widehat D^\ast P_{r\B}(\ox)(\ox)=\emptyset.$
\end{Corollary}
	{\bf Proof.}
(a) Using \eqref{kq} with $y=0$, we have  
$$\begin{array}{rl}
\widehat D^\ast P_{r\B}(\ox)(\theta)&=\{z\in H: a(\theta)\leq a(z)\leq 0, \  o(z)=o(\theta)\}\\
&=\{z\in H:  a(z)= 0, \  o(z)=\theta\}\\
&=\{\theta\}.
\end{array}$$

(b)  By \eqref{kq}, $z\in \widehat D^\ast P_{r\B}(\ox)(y) \Longleftrightarrow a(y)\leq a(z)\leq 0, \ o(z)=o(y).$ So,
$$\begin{array}{rl}\theta\in \widehat D^\ast P_{r\B}(\ox)(y) &\Longleftrightarrow a(y)\leq a(\theta)\leq 0, \ o(\theta)=o(y)\\
	&\Longleftrightarrow a(y)\leq 0, \ o(y) =\theta\\
	&\Longleftrightarrow y=a(y)\ox, \ \mbox{with}\ \ a(y)=\dfrac{\la y,\ox\ra}{\|\ox\|^2}\leq 0.
\end{array}	.$$

(c) Using \eqref{kq} with $y=\ox$, we have  
$$\begin{array}{rl}
	\widehat D^\ast P_{r\B}(\ox)(\ox)&=\left\{z\in H: \la \ox,\ox\ra\leq \la z,\ox\ra \leq 0,\ \    \ox-z=\dfrac{\la \ox-z,\ox\ra}{\|\ox\|^2}\ox\right\}\\
	&=\left\{z\in H:  1\leq \la z,\ox\ra\leq 0, \  z=\dfrac{\la z,\ox\ra}{\|\ox\|^2}\ox\right\}\\
	&=\emptyset.
\end{array}$$\hfill $\square$\\

To proceed, we need the following result, which provides the relationship between $P_{r\B}$ and $P_{r\B}(x-c)$.
\begin{Lemma} \label{Lem33}
	Let $C=\B(c,r)$, where $c\in H, r>0$. Then, we have
	$$P_C(x)=P_{r\B}(x-c)+c,\ \ \ \mbox{for all}\ x\in H.$$
\end{Lemma}

We now arrive at the main result of this section, the following theorem provides the formula for calculating the regular coderivatives of the metric projection onto the closed balls with center at arbitrarily given point $c$ in $H$.
\begin{Theorem} \label{Thm34}
Let $H$ be a Hilbert space. For any  $c\in H, r>0$ put $C=\B(c,r)$. Then, the regular coderivative of $P_{C}:H\rightarrow C$   at a point $\ox\in H$ is given below.

(i)  If $\ox\in C^o$, $$\widehat D^\ast P_{C}(\ox)(y)=\{y\} \ \ \mbox{for every}\ y\in H.$$

(ii) If $\ox\in H\backslash C$, then
$$\widehat D^\ast P_{C}(\ox)(y)=\left\{\dfrac{r}{\|\ox-c\|}\Big(y-\dfrac{\la \ox-c, y\ra}{\|\ox-c\|^2}(\ox-c)\Big)\right\}\ \ \mbox{for every} \ y\in H.$$

(iii) If $\ox\in \mbox{\rm bd}C$, then
$$\widehat D^\ast P_{C}(\ox)(y)=\left\{z\in H: \la y,\ox-c\ra\leq \la z,\ox-c\ra \leq 0, \  y-z=\dfrac{\la y-z,\ox-c\ra}{\|\ox-c\|^2}(\ox-c)\right\}.$$
\end{Theorem}

{\bf Proof.} By Lemma \ref{Lem33}, we have $$P_C=P_{r\B}\circ\psi + h$$
where $\psi,\ph,h:H\to H$ with $\psi(x)=x-c, \ph(x)=x+c$, and $h(x)=c$ for all $x\in H.$ Note that $\psi,\ph,h$ are strictly differentiable at any $x\in H$ and $\nabla\psi(x)=\nabla\ph(x)=I_H, \ \nabla h(x)=\theta_H.$\\

(i). For any given $\ox\in C^o$, we have $\ox-c\in r\B^o$, by Theorem \ref{Thm1} and Propositions \ref{Pro13}, \ref{Pro23}, \ref{Pro24}, we get 
$$\begin{array} {rl}\widehat D^\ast P_C(\ox)(y)&= \nabla\psi(\ox)\widehat D^\ast P_{r\B}(\psi(\ox))(y)+\nabla h(\ox)(y)\\
&=\widehat D^\ast P_{r\B}(\ox-c)(y)+\theta_H(y)\\
&=\{y\}\ \ \ \ \ \ \ \ \ \ \ \ \ \ \ \ \ \ \ \ \ \ \ \ \ \ \ \ \ \ \ \ \  \ \ \ \ \ \ \  \mbox{for every }\ y\in H. \end{array}$$

(ii).  For $\ox\in H\backslash C$, we have $\ox-c\in H\backslash r\B$, by Theorem \ref{Thm1} and Propositions \ref{Pro13}, \ref{Pro23}, \ref{Pro24}, we get 
$$\begin{array} {rl}\widehat D^\ast P_C(\ox)(y)&= \nabla\psi(\ox)\widehat D^\ast P_{r\B}(\psi(\ox))(y)+\nabla h(\ox)(y)\\
&=\widehat D^\ast P_{r\B}(\ox-c)(y)+\theta_H(y)\\
	&=\left\{\dfrac{r}{\|\ox-c\|}\Big(y-\dfrac{\la \ox-c, y\ra}{\|\ox-c\|^2}(\ox-c)\Big)\right\}+\{\theta\}\\
	&=\left\{\dfrac{r}{\|\ox-c\|}\Big(y-\dfrac{\la \ox-c, y\ra}{\|\ox-c\|^2}(\ox-c)\Big)\right\},\ \ \ \ \mbox{for all} \ y\in H. \end{array}$$

(iii). For $\ox\in \mbox{bd}C$, we have $\ox-c\in \mbox{bd}r\B$, by Theorem \ref{Thm1} and Propositions \ref{Pro13}, \ref{Pro23}, \ref{Pro24}, we get 
$$\begin{array} {rl}\widehat D^\ast P_C(\ox)(y)&= \nabla\psi(\ox)\widehat D^\ast P_{r\B}(\psi(\ox))(y)+\nabla h(\ox)(y)\\
	&=\widehat D^\ast P_{r\B}(\ox-c)(y)+\theta_H(y)\\
	&=\left\{y\in H: \la y,\ox-c\ra\leq \la z,\ox-c\ra \leq 0, \  y-z=\dfrac{\la y-z,\ox-c\ra}{\|\ox-c\|^2}(\ox-c)\right\}+\{\theta\}\\
	&=\left\{y\in H: \la y,\ox-c\ra\leq \la z,\ox-c\ra \leq 0, \  y-z=\dfrac{\la y-z,\ox-c\ra}{\|\ox-c\|^2}(\ox-c)\right\}, \end{array}$$ $\mbox{for all} \ y\in H$ \hfill $\square$

To end this section, we give some examples to demonstrate the results of Theorem \ref{Thm34}.
\begin{Example} {\rm 
	Let $H=\R$, for any $c\in \R, r>0, \ C=\B(c,r)=[c-r;c+r],$ $C^o=(c-r;c+r)$ and $\mbox{bd}C=\{c-r;c+r\}$. For $x\in \R$, we have
	$$\widehat D^\ast P_C(x)(y)=\begin{cases}
		\{y\}\ \ \ \ \ \ \ \ \ \ \ \mbox{if}\ \ x\in (c-r; c+r),\ y\in\R,\\
		\{0\}  \ \ \ \ \ \ \ \ \  \ \ \mbox{if}\ \ x\not\in [c-r; c+r], \ y\in\R,\\
		[y; 0]\ \ \ \ \ \ \ \ \  \mbox{if}\ \ x= c+r, \ y\leq 0,\\
		[0; y]\ \ \ \ \ \ \ \ \  \mbox{if}\ \ x= c-r, \ y\geq 0,\\
		\emptyset \ \ \ \ \ \ \ \ \ \  \ \  \ \ \mbox{if}\ \ x=c-r, y<0 \ \mbox{or}\ x= c+r, \ y> 0.
	\end{cases}$$
}
\end{Example}

\begin{Example} {\rm
		Let $H=\R^2$ and $C=\B(c,r)$ is the closed ball in $\R^2$ with radius $r>0$ and center $c=(c_1,c_2)\in\R^2$. For any $x=(x_1,x_2)\in \R^2,$ we have
	$$\begin{array} {rl} &\widehat D^\ast P_C(x)(y)\\
		&=\begin{cases}
		\left\{y\right\}\ \ \ \ \ \ \ \ \ \ \ \ \ \ \ \ \ \ \ \ \ \ \ \ \ \ \ \ \ \ \ \ \ \ \ \ \ \ \ \ \ \ \ \ \ \ \ \ \ \ \ \  \ \ \ \ \ \ \ \ \mbox{if}\ \ (x_1-c_1)^2+(x_2-c_2)^2<r^2\\
		\left\{\dfrac{r}{\sqrt{(x_1-c_1)^2+(x_2-c_2)^2}}\bigg((y_1,y_2)-\dfrac{y_1(x_1-c_1)+y_2(x_2-c_2)}{(x_1-c_1)^2+(x_2-c_2)^2}(x_1-c_1,x_2-c_2) \bigg)\right\}\\  \ \ \ \ \ \ \ \ \ \ \ \ \ \ \ \ \ \ \ \ \  \ \ \ \ \ \ \ \ \ \ \ \ \ \ \ \ \ \ \ \ \ \ \ \ \ \ \ \ \ \ \ \ \ \ \ \ \ \ \  \ \ \mbox{if}\ \ (x_1-c_1)^2+(x_2-c_2)^2>r^2\\
		\left\{\begin{array} {rl}&(z_1,z_2)\in\R^2\big|y_1(x_1-c_1)+y_2(x_2-c_2)\leq z_1(x_1-c_1)+z_2(x_2-c_2)\leq 0,\\
			&(y_1-z_1, y_2-z_2)=\dfrac{(y_1-z_1)(x_1-c_1)+(y_2-z_2)(x_2-c_2)}{(x_1-c_1)^2+(x_2-c_2)^2}\big(x_1-c_1, x_2-c_2\big)\end{array} \right\}\\ \ \ \ \ \ \ \ \ \ \ \ \ \ \ \ \ \ \ \ \ \  \ \ \ \ \ \ \ \ \ \ \ \ \ \ \ \ \ \ \ \ \ \ \ \ \ \ \ \ \ \ \ \ \ \ \ \ \ \ \  \ \  \mbox{if}\ \ (x_1-c_1)^2+(x_2-c_2)^2=r^2
	\end{cases}\\&  \mbox{for every}\ y=(y_1,y_2)\in\R^2.\end{array}$$
}
\end{Example}


\section{Graphical derivative of the metric projection operator onto ball in Hilbert spaces}\label{Sec3}
\setcounter{equation}{0}

\begin{Definition} {\rm (see \cite{M18,RW98})}
	Let $\Omega$ be a nonempty subset of  a Banach  space $X$. The (Bouligand-Severi) tangent/contingent cone to the set $\Omega$ at $\ox\in\Omega$ is defined by
	$$T_\Omega(\ox):=\big\{z\in X|\ \mbox{there exist}\ t_k\downarrow 0, z_k\rightarrow z \ \mbox{with}\ \ox + t_kz_k \in \Omega \ \mbox{for all}\ k\in\mathbb{N}\big\}.$$
\end{Definition}
If $\ox\not\in\Omega$, one puts
$T_\Omega(\ox)=\emptyset$ by convention. When the set $\Omega$ is convex, the 	above tangent cone and normal cone reduce to the tangent cone and normal cone inthe sense of classical convex analysis.

\begin{Definition}{\rm (see \cite{M18,RW98})}
	Let $X, Y$ be Banach spaces and $\Phi: X\rightrightarrows Y$ be a set-valued mapping with its graph {\rm gph}$\Phi: =\{(x,y)|y\in\Phi(x)\}$ and its domain {\rm Dom}$\Phi:=\{x|\Phi(x)\not=\emptyset\}$. Given a point $\ox\in$ {\rm Dom}$\Phi$, the graphical derivative of $\Phi$ at $\ox$ for $\oy\in\Phi(\ox)$ is	the set-valued mapping $D\Phi(\ox|\oy):X \rightrightarrows Y$ defined by
	$$D\Phi(\ox|\oy)(v):=\big\{w\in X|(v,w)\in T_{\mbox{\rm gph}\Phi}(\ox, \oy)\big\} \ \mbox{for all}\ v\in H,$$
	that is, {\rm gph}$D\Phi(\ox|\oy):=T_{\mbox{\rm gph}\Phi}(\ox, \oy)$
	
	In the case $\Phi(\ox)=\{\oy\},$ one writes $D\Phi(\ox)$ for $D\Phi(\ox|\oy)$.
\end{Definition}

\begin{Theorem} \label{Thm222}
 Let $X, Y$ be Banach spaces, $f:X\to Y$ be a single-valued mapping. Suppose that $f$ is Lipschitz continuous around $\ox\in X$ and Gâteaux directional differentiable at point $\ox$ along direction $u\in X.$ Then, $Df(\ox)(u)=\{f'(\ox)(u)\}.$
\end{Theorem}
{\bf Proof.} Suppose that $v\in Df(\ox)(u)$. Then, there exists $t_k\downarrow 0, (u_k, v_k)\rightarrow (u,v)$ with $(\ox,f(\ox))+t_k(u_k, v_k)\in \mbox{\rm gph}f$ for all $k\in \mathbb{N}.$ For every  $k\in \mathbb{N},$ we have  $f(\ox)+t_kv_k=f(\ox+t_ku_k)$, it means $$v_k=\dfrac{f(\ox+t_ku_k)-f(\ox)}{t_k}.$$
So, \begin{equation} \label{41}\begin{array}{rl}
	v&=\lim\limits_{k\rightarrow \infty} v_k\\
	&=\lim\limits_{k\rightarrow \infty} \dfrac{f(\ox+t_ku_k)-f(\ox)}{t_k}\\
		&=\lim\limits_{k\rightarrow \infty} \dfrac{f(\ox+t_ku)-f(\ox)}{t_k} + \lim\limits_{k\rightarrow \infty} \dfrac{f(\ox+t_ku_k)-f(\ox+t_ku)}{t_k}.
\end{array}\end{equation}
By the Gâteaux direction differentiability of $f$, we have
\begin{equation} \label{42}\lim\limits_{k\rightarrow \infty} \dfrac{f(\ox+t_ku)-f(\ox)}{t_k}=\lim\limits_{t\rightarrow 0} \dfrac{f(\ox+tu)-f(\ox)}{t}=f'(\ox)(u).\end{equation}
Thanks to the Lipschitz continuity of $f$ around $\ox$, there exists $L\geq 0$ such that $$\|f(\ox+t_ku_k)-f(\ox+t_ku)\|\leq Lt_k\|u_k-u\|.$$
Therefore, \begin{equation} \label{43}
0\leq \left\|\dfrac{f(\ox+t_ku_k)-f(\ox+t_ku)}{t_k}\right\|\leq \dfrac{Lt_k\|u_k-u\|}{t_k} = L\|u_k-u\| \rightarrow 0\ \ \ \  \mbox{as} \ k\rightarrow \infty.
\end{equation}
Combining \eqref{41} with \eqref{42} and \eqref{43}, we obtain $v=f'(\ox)(u).$

Next, we prove the converse, assuming $v=f'(\ox)(u)$, we have
$$v=\lim\limits_{t\rightarrow 0} \dfrac{f(\ox+tu)-f(\ox)}{t}.$$
Choose $t_k\rightarrow 0, u_k=u$ and $v_k=\dfrac{f(\ox+t_ku_k)-f(\ox)}{t_k}$ for all $k\in\mathbb{N}$.\\
 We get $t_k\downarrow 0, (u_k, v_k)\rightarrow (u,v)$ and $f(\ox)+t_kv_k=f(\ox+t_ku_k)$ for all $k\in \mathbb{N},$ it means $(\ox, f(\ox))+t_k(u_k, v_k)\in \ \mbox{\rm gph}f.$ This implies that $v\in Df(\ox)(u)$.\hfill $\square$

\begin{Proposition} \label{Pro111}{\rm (see \cite{LLX23})}
Let $H$ be a Hilbert space. For any  $c\in H, r>0$ put $C=\B(c,r)$. Then, $P_C$ is Gâteaux directionally differentiable on $H$ such that, for every $u\in H$, $P_C$ has the following Gâteaux directional differentiability.

(i) For any $\ox\in C^o$, we have  $$ P'_{C}(\ox)(u)=u.$$

(ii) For any $\ox\in H\backslash C$, we have  
$$P'_{C}(\ox)(u)=\dfrac{r}{\|\ox-c\|}\Big(u-\dfrac{\la \ox-c, u\ra}{\|\ox-c\|^2}(\ox-c)\Big).$$

(iii) For any $\ox\in \mbox{\rm bd}C$, we have
$$P^{'}_{C}(\ox)(u)=\begin{cases}
	u-\dfrac{1}{r^2}\la \ox-c,u\ra(\ox-c) \ \ \ \mbox{if} \ \ \ u\in\ox^\uparrow_{(c,r)};\\u \ \ \ \  \ \ \ \ \ \ \ \ \ \ \ \ \ \ \ \ \ \ \ \ \ \  \ \ \ \  \mbox{if} \  \ \  u\in\ox^\downarrow_{(c,r)}.
\end{cases}$$  
\end{Proposition}

\begin{Theorem} \label{kqc}
Let $H$ be a Hilbert space. For any  $c\in H, r>0$ put $C=\B(c,r)$. Then, the graphical derivative of $P_C:H\rightarrow C$   at a point $\ox\in H$ is given below.

(i)  If $\ox\in C^o$, $$D P_{C}(\ox)(y)=\{y\} \ \ \mbox{for every}\ y\in H.$$

(ii) If $\ox\in H\backslash C$, then
$$D P_{C}(\ox)(y)=\left\{\dfrac{r}{\|\ox-c\|}\Big(y-\dfrac{\la \ox-c, y\ra}{\|\ox-c\|^2}(\ox-c)\Big)\right\}\ \ \mbox{for all} \ y\in H.$$

(iii) If $x\in\mbox{\rm bd}C$, we have
$$D P_{C}(\ox)(y)=\begin{cases}
	\left\{y-\dfrac{1}{r^2}\la \ox-c,y\ra(\ox-c)\right\} \ \ \ \mbox{if} \ \ \ y\in\ox^\uparrow_{(c,r)};\\\left\{y\right\}\ \ \ \  \ \ \ \ \ \ \ \ \ \ \ \ \ \ \ \ \ \ \ \ \ \  \ \ \ \  \mbox{if} \  \ \ \ y\in\ox^\downarrow_{(c,r)}.
\end{cases}$$  
\end{Theorem}
{\bf Proof.} The assertions of the Theorem are derived directly from Propositions \ref{Pro111} and Theorem \ref{Thm222}, with note that the function $P_C$ is Lipschitz continuous  on $H$.\hfill $\square$

\begin{Example} {\rm 
		Let $H=\R$, for any $c\in \R, r>0, \ C=\B(c,r)=[c-r;c+r],$ $C^o=(c-r;c+r)$ and $\mbox{bd}C=\{c-r;c+r\}$. For $x=c+r,\ x^\uparrow_{(c,r)}=\{y\in\R| y\geq 0\},\ x^\downarrow_{(c,r)}=\{y\in\R| y< 0\}$ and $x=c-r,\ x^\uparrow_{(c,r)}=\{y\in\R| y\leq 0\},\ x^\downarrow_{(c,r)}=\{y\in\R| y> 0\}$. So, by Theorem \ref{kqc}, we have
		$$DP_C(x)(y)=\begin{cases}
			\{y\}\ \ \ \mbox{if}\  \ x\in (c-r; c+r),\ y\in\R\ \mbox{or}\ x=c+r, y< 0 \ \mbox{or} \ x=c=r, y>0,\\
			\{0\} \ \ \  \mbox{if}\ \ x\not\in [c-r; c+r], \ y\in\R \ \mbox{or}\ x=c+r, y\geq0\ \mbox{or}\ x=c-r, y\leq 0.
		\end{cases}$$
	}
\end{Example}

\begin{Example} {\rm
		Let $H=\R^2$ and $C=\B(c,r)$ is the closed ball in $\R^2$ with radius $r>0$ and center $c=(c_1,c_2)\in\R^2$. For any $x=(x_1,x_2)\in \mbox{\rm bd}C,$ we have $x^\uparrow_{(c,r)}=\{y\in\R^2|(x_1-c_1)y_1+(x_2-c_2)y_2\geq 0\},\ \ x^\downarrow_{(c,r)}=\{y\in\R^2|(x_1-c_1)y_1+(x_2-c_2)y_2< 0\}.$ For $y=(y_1,y_2)\in\R^2$, by Theorem \ref{kqc}, we obtain 
		$$\begin{array} {rl} &DP_C(x)(y)\\
			&=\begin{cases}
				\{y\} \ \ \ \ \ \ \ \ \ \ \ \  \ \ \ \ \ \ \ \begin{array}{rl} &\mbox{if}\  (x_1-c_1)^2+(x_2-c_2)^2<r^2, y=(y_1,y_2)\in \R^2\\ 
					&\mbox{or}\ (x_1-c_1)^2+(x_2-c_2)^2=r^2,\ (x_1-c_1)y_1+(x_2-c_2)y_2\geq 0;\end{array}\\
				\left\{\dfrac{r}{\sqrt{(x_1-c_1)^2+(x_2-c_2)^2}}\bigg((y_1,y_2)-\dfrac{y_1(x_1-c_1)+y_2(x_2-c_2)}{(x_1-c_1)^2+(x_2-c_2)^2}(x_1-c_1,x_2-c_2) \bigg)\right\}\\  \ \ \ \ \ \ \  \ \ \ \ \ \  \ \ \ \ \ \ \ \ \ \  \begin{array}{rl}&\mbox{if}\  (x_1-c_1)^2+(x_2-c_2)^2>r^2,\  y=(y_1,y_2)\in \R^2\\
					&\mbox{or}\ \ (x_1-c_1)^2+(x_2-c_2)^2=r^2,\ (x_1-c_1)y_1+(x_2-c_2)y_2\geq 0. \end{array}
				
			\end{cases}
		\end{array}$$
	}
\end{Example}

\section{Concluding remarks}

The main result of this paper exhibits formula for computing the regular derivative of the metric projection operator onto closed balls centered at origin, and we generalize this result for the metric projection operator onto closed balls with center at arbitrarily given point $c$ in Hilbert spaces. In addition, we also prove the formula for calculating the graphical derivative of the metric projection operator in this case. In the near future, we will continue to research other generalized differential structures of metric projection operators in Hilbert spaces.

.
\section*{Acknowledgments}
The author would like to thank Vietnam Institute for Advanced Study in Mathematics for hospitality during his post-doctoral fellowship of the Institute in 2022--2023.

\end{document}